\documentclass[12pt,oneside,reqno]{amsart}
\usepackage{graphicx}
\usepackage{mathrsfs}
\usepackage{stmaryrd}
\usepackage{amsfonts}
\usepackage{enumerate,amsmath,amssymb,amsthm}
\newcommand{\dif}{\mathrm{d}}

\newcommand{\be}{\begin{eqnarray}}
\newcommand{\ee}{\end{eqnarray}}
\newcommand{\ce}{\begin{eqnarray*}}
\newcommand{\de}{\end{eqnarray*}}
\newtheorem{theorem}{Theorem}[section]
\newtheorem{lemma}[theorem]{Lemma}
\newtheorem{remark}[theorem]{Remark}
\newtheorem{definition}[theorem]{Definition}
\newtheorem{proposition}[theorem]{Proposition}
\newtheorem{Examples}[theorem]{Examples}
\newtheorem{corollary}[theorem]{Corollary}

\def\[{{\Big[}}
\def\]{{\Big]}}
\def\<{{\langle}}
\def\>{{\rangle}}
\def\({{\Big(}}
\def\){{\Big)}}

\def\no{\nonumber}
\def\bt{\begin{theorem}}
\def\et{\end{theorem}}
\def\bl{\begin{lemma}}
\def\el{\end{lemma}}
\def\br{\begin{remark}}
\def\er{\end{remark}}
\def\bx{\begin{Examples}}
\def\ex{\end{Examples}}
\def\bd{\begin{definition}}
\def\ed{\end{definition}}
\def\bp{\begin{proposition}}
\def\ep{\end{proposition}}
\def\bc{\begin{corollary}}
\def\ec{\end{corollary}}

\def\cM{{\mathcal M}}

\def\cS{{\mathcal S}}

\def\cW{{\mathcal W}}

\def\mE{{\mathbb E}}

\def\mN{{\mathbb N}}

\def\mP{{\mathbb P}}

\def\mR{{\mathbb R}}

\def\sB{{\mathscr B}}

\def\sF{{\mathscr F}}

\def\sW{{\mathscr W}}

\def\geq{\geqslant}
\def\leq{\leqslant}

\begin{document}

\allowdisplaybreaks

\title{Euler-Maruyama Approximations for Stochastic McKean-Vlasov Equations with Non-Lipschitz Coefficients$^{\bigstar}$}

\author{Xiaojie Ding$^{\sharp}$ and Huijie Qiao$^{*}$}



\thanks{${\bigstar}$ This work was supported by NSF of China (No. 11001051, 11371352) and China Scholarship Council under Grant No. 201906095034.}

\thanks{${\sharp}$ 220171469@seu.edu.cn}

\thanks{${*}$ Corresponding author: hjqiaogean@seu.edu.cn}


\date{}
\dedicatory{School of Mathematics,
Southeast University,\\
Nanjing, Jiangsu 211189, P.R.China}

\begin{abstract}
In this paper we study a type of stochastic McKean-Vlasov equations with non-Lipschitz coefficients. Firstly, by an Euler-Maruyama approximation the existence of its weak solutions is proved. And then we observe the pathwise uniqueness of its weak solutions. Finally, it is shown that the Euler-Maruyama approximation has an optimal strong convergence rate.
\end{abstract}

\maketitle \rm
\section{Introduction}
Given $T>0$. Suppose that a complete filtered probability space
$(\Omega,\mathscr{F},\{\mathscr{F}_t\}_{t\in[0,T]},\mP)$ and a standard $d$-dimensional Brownian motion $W_t$ on the probability space are given. Consider the following stochastic McKean-Vlasov equation(SMVE) on $\mR^d$:
\be\left\{\begin{array}{ll}
X_t=\xi+\int_0^tb(X_s,\mu_s)ds+\int_0^t\sigma(X_s,\mu_s)d W_s,\\
\mu_s= $probability  distribution of\quad$X_s,
\end{array}
\label{eq1}
\right.
\ee
where $\xi$ is a $\sF_0$-measurable random variable, and the coefficients $b:\mR^d\times\cM_{\lambda^2}(\mR^d)\rightarrow{\mR^d}, \sigma:\mR^d\times\cM_{\lambda^2}(\mR^d)\rightarrow{\mR^d}\times{\mR^d}$ are Borel measurable ($\cM_{\lambda^2}(\mR^d)$ is defined in Section \ref{nn}).

If $b$ and $\sigma$ in Eq.$(\ref{eq1})$ are independent of the probability distribution $\mu_t$ of the process at time $t$, Eq.$(\ref{eq1})$ is a standard Markov process and has been well studied in the literature (c.f.\cite{hzy}\cite{dz}). Besides, there are stochastic differential equations whose coefficients depend not only on the process but also on the probability  distribution of the process at time $t$ as indicated in Eq.$(\ref{eq1})$. The study on SMVEs was initiated by Henry P. McKean \cite{MC} who was inspired by Kac's program in Kinetic Theory. And then there have been numerous results (c.f.\cite{Sz}). Let us mention some works. Recently,  Huang-Wang \cite{hw} studied the existence and the uniqueness of strong solutions for Eq.$(\ref{eq1})$ under some integrable conditions. Besides, if the diffusion coefficient is independent of $\mu_t$, the second named author \cite{q1} showed that under Lipschitz and linear growth conditions, Eq.$(\ref{eq1})$ has a unique mild solution in a real separated Hilbert space, and the Euler approximation of the mild solution converges to itself. Later, under more general conditions than that in \cite{q1}, Govindan-Ahmed \cite{GA} proved Eq.$(\ref{eq1})$ has a unique mild solution, and the Yosida appromiximation of the mild solution converges to itself. If $b$ and $\sigma$ depend on $\mu_t$ as follows:
\ce
&&\int_0^tb[X_s,\mu_s]ds=\int_0^t\int_{\mR^d}b(X_s,y)\mu_s(dy)ds,\\
&&\int_0^t\sigma[X_s,\mu_s]ds=\int_0^t\int_{\mR^d}\sigma(X_s,y)\mu_s(dy)ds,
\de
where $b:\mR^d\times\mR^d\rightarrow{\mR^d}, \sigma:\mR^d\times\mR^d\rightarrow{\mR^d}\times{\mR^d}$ are Borel measurable, Sznitman  \cite{Sz} investigated the existence and the uniqueness of strong solutions for Eq.$(\ref{eq1})$ with a fixed point argument if the coefficients are globally Lipschitz continuous. Recently, Chi \cite{chm} proved that if the coefficients are continuous and satisfy linear growth condition, a weak solution of the multivalued SMVE exists by the Euler approximation. 

In this paper, we study Eq.$(\ref{eq1})$ under non-Lipschitz conditions. Firstly, we establish the weak existence of Eq.$(\ref{eq1})$ under a linear growth condition. Next, the pathwise uniqueness is obtained under two non-Lipschitz conditions. Thus, by the weak existence and the pathwise uniqueness, we prove the existence and the uniqueness of a strong solution for Eq.$(\ref{eq1})$. Finally, the convergence rate of the Euler-Maruyama approximation is discussed.

It is worthwhile to mentioning our conditions and methods. We give two non-Lipschitz conditions which can not be covered by the conditions in \cite{hw}. Moreover, our conditions are more straight than that in \cite{hw}. Besides, we prove the existence of martingale solutions of Eq.$(\ref{eq1})$ by an Euler-Maruyama approximation , which implies its weak existence. Thus, a number of complex calculation, as that in \cite{q3} \cite{zx}, is avoided.

The rest of the paper is organized as follows. In Section \ref{fram}, we recall some basic notations, and give some necessary concepts and assumptions. And then we prove the existence and the uniqueness of a strong solution of Eq.$(\ref{eq1})$ in Section \ref{exun}. In Section \ref{cora}, the convergence rate of the Euler-Maruyama approximation is investigated.

The following convention will be used throughout the paper: $C$ with or without indices will denote different positive constants whose values may change from one place to
another.

\section{The Framework}\label{fram}

In the section, we recall some basic notations, and give some necessary concepts and assumptions.

\subsection{Notations}\label{nn}

In the subsection, we introduce notations used in the sequel.

Let $C(\mR^d)$ be the space of continuous functions on $\mR^d$. And let $C_0^k(\mR^d)$ be the collection of all continuous functions which have bounded, continuous partial derivatives of every order up to $k$ where $k$ is a positive integer. Let $\partial_{ij}$ denote the differentiation with respect to the coordinates with corresponding numbers (e.g. $\partial_{ij}(f):=\frac{\partial^2f(x)}{\partial{x^i}\partial{x^j}}$). Let $\sB(\mR^d)$ be the Borel $\sigma$-algebra on $\mR^d$ and $\cM({\mR^d})$ be the space of all probability measures defined on $\sB(\mR^d)$ carrying the usual topology of weak convergence.

For convenience, we shall use $\mid\cdot\mid$
and $\parallel\cdot\parallel$  for norms of vectors and matrices, respectively. Furthermore, let $\langle\cdot$ , $\cdot\rangle$ denote the scalar product in $\mR^d$. Let $A^*$ denote the transpose of the matrix $A$.

Define the Banach space
$$
C_\rho(\mR^d):=\left\{{\varphi\in{C(\mR^d)},\parallel{\varphi}\parallel_{C_\rho(\mR^d)}
=\sup_{x\in{\mR^d}}\frac{\mid{\varphi(x)}\mid}{(1+\mid{x}\mid)^2}+\sup_{x\neq{y}}\frac{\mid{\varphi(x)-\varphi(y)}\mid}{\mid{x-y}\mid}<\infty}\right\}.
$$
Let $\cM_{\lambda^2}^s(\mR^d)$ be the Banach space of signed measures $m$
on $\sB(\mR^d)$ satisfying
\ce
\|m\|_{\lambda^2}^2:=\int_{\mR^d}(1+\mid{x}\mid)^2\,|m|(\dif x)<\infty,
\de
where $|m|=m^{+}+m^{-}$ and $m=m^{+}-m^{-}$ is the Jordan decomposition of $m$. Let
$\cM_{\lambda^2}(\mR^d)=\cM_{\lambda^2}^s(\mR^d)\bigcap\cM(\mR^d)$ be the set of probability
measures on $\sB(\mR^d)$. We put on $\cM_{\lambda^2}(\mR^d)$ a topology induced by the
following metric:
\ce
\rho(\mu,\nu):=\sup_{\parallel{\varphi}\parallel_{C_\rho(\mR^d)\leq1}}\left|{\int_{\mR^d}\varphi(x)\mu(dx)-\int_{\mR^d}\varphi(x)\nu(dx)}\right|.
 \de
Then $(\cM_{\lambda^2}(\mR^d),\rho)$ is a complete metric space.

\subsection{Some concepts}\label{conc}

In the subsection, we introduce the concepts of strong solutions, weak solutions and pathwise uniqueness. Consider Eq.$(\ref{eq1})$, i.e.
\ce\left\{\begin{array}{ll}
X_t=\xi+\int_0^tb(X_s,\mu_s)ds+\int_0^t\sigma(X_s,\mu_s)d W_s,\\
\mu_s= $probability  distribution of\quad$X_s.
\end{array}
\right.
\de

\bd\label{sd3}
We say that Eq.$(\ref{eq1})$ admits a strong solution with the initial value $\xi$ if there exists a continuous process
$X=\{X_t;0\leq{t}\leq{T}\}$ on $(\Omega,\mathscr{F},\{\mathscr{F}_t\}_{t\in[0,T]},\mP)$ such that

(i) $\mP(X_0=\xi)=1$,

(ii) $X_t\in{\mathscr{F}_t^W}$, where $\{\mathscr{F}_t^W\}_{t\in[0,T]}$ stands for the $\sigma$-field filter generated by $W$,

(iii) it holds that
\ce
\int_0^t(\mid{b(X_s,\mu_s)}\mid+\parallel{\sigma(X_s,\mu_s)}\parallel^2)ds<+\infty,\quad a.s.{\mP},
\de
and
\ce
X_t=\xi+\int_0^tb(X_s,\mu_s)ds+\int_0^t\sigma(X_s,\mu_s)d W_s, \quad 0\leq{t}\leq{T}. \de
\ed

From the above definition, we know that $\mu_0=\mP\circ\xi^{-1}$.

\bd\label{sd2}
We say that Eq.$(\ref{eq1})$ admits a weak solution with the initial law $\mu_0$ if there exists a stochastic space $\hat{\cS}:=(\hat{\Omega},\hat{\mathscr{F}},\{\hat{\mathscr{F}_t\}}_{t\in[0,T]},\hat{\mP})$, a d-dimensional standard Brownian motion $\hat{W}$ as well as a $\{\hat{\mathscr{F}_t\}}_{t\in[0,T]}$-adapted process $\hat{X}$ defined on $\hat{\cS}$ such that

(i) $\hat{\mP}\circ\hat{X}^{-1}_0=\mu_0$,

(ii) it holds that
\ce
\int_0^t(\mid{b(\hat{X}_s,\hat{\mu}_s)}\mid+\parallel{\sigma(\hat{X}_s,\hat{\mu}_s)}\parallel^2)ds<+\infty,\quad a.s.\hat{\mP},
\de
and
\ce
\hat{X}_t=\hat{X}_0+\int_0^tb(\hat{X}_s,\hat{\mu}_s)ds+\int_0^t\sigma(\hat{X}_s,\hat{\mu}_s)d \hat{W}_s,  \quad 0\leq{t}\leq{T}.
\de
\ed

Such a weak solution  is denoted by ($\hat{\cS}; \hat{W}, \hat{X}$).

\bd\label{sd4}
(Pathwise Uniqueness)\quad Suppose ($\hat{\cS}; \hat{W}, \hat{X}^1$) and ($\hat{\cS}; \hat{W}, \hat{X}^2$) are two weak solutions with $\hat{X}^1_0=\hat{X}^2_0$. If $\hat{\mP}(\hat{X}^1_t=\hat{X}^2_t,t\geq0)=1$, then we say that the pathwise uniqueness holds for Eq.$(\ref{eq1})$.
\ed

\subsection{Some assumptions}\label{ass}

In the subsection, we give out some assumptions:

\begin{enumerate}[($\bf{H}_1$)]
\item The functions $b, \sigma$ are continuous in $(x,\mu)$ and satisfy for $(x,\mu)\in\mR^{d}\times{\cM_{\lambda^2}(\mR^d)}$
\be
{\mid{b(x,\mu)}\mid}^2+{\parallel{\sigma(x,\mu)}\parallel}^2\leq{L_1(1+\mid{x}\mid^2+\parallel{\mu}\parallel^2_{\lambda^2})}, 
\label{eq3}
\ee
where $L_1>0$ is a constant.
\end{enumerate}
\begin{enumerate}[($\bf{H}_2$)]
\item The functions $b,\sigma$ satisfy for $(x_1,\mu_1), (x_2,\mu_2)\in\mR^{d}\times{\cM_{\lambda^2}(\mR^d)}$
\ce
2\langle{x_1-x_2,b(x_1,\mu_1)-b(x_2,\mu_2)}\rangle+\parallel{\sigma(x_1,\mu_1)-\sigma(x_2,\mu_2)}\parallel^2\leq{L_2\(\kappa_1(|x_1-x_2|^2)+\kappa_1\left(\rho^2(\mu_1,\mu_2)\right)\)},
\de
where $L_2>0$ is a constant, and $\kappa_i(x), i=1, 2$ are two positive, strictly increasing, continuous concave function and satisfies $\kappa_i(0)=0$,  $\int_{0^+}\frac{1}{\kappa_1(x)+\kappa_2(x)}dx=\infty$.
\end{enumerate}
\begin{enumerate}[($\bf{H}_2'$)]
\item The functions $b$ and $\sigma$ satisfy for $(x_1,\mu_1), (x_2,\mu_2)\in\mR^{d}\times{\cM_{\lambda^2}(\mR^d)}$
\ce
\mid b(x_1,\mu_1)-b(x_2,\mu_2)\mid&\leq&\lambda_1\(|x_1-x_2|\gamma_1(|x_1-x_2|)+\rho(\mu_1,\mu_2)\),\\
\parallel{\sigma(x_1,\mu_1)-\sigma(x_2,\mu_2)}\parallel^2&\leq&\lambda_2\(|x_1-x_2|^2\gamma_2(|x_1-x_2|)+\rho^2(\mu_1,\mu_2)\),
\de
where $\lambda_i>0$ is a constant and $\gamma_i(x)$ is a positive continuous function, bounded on $[1,\infty)$ and satisfying
\ce
\lim_{x\downarrow0}\frac{\gamma_i(x)}{\log(x^{-1})}=\delta_i<\infty,\quad\quad i=1,2.
\de
\end{enumerate}

\br\label{h2h2}
If $b(x,\mu)$ satisfies ($\bf{H}_2'$), it holds that for $(x_1,\mu_1), (x_2,\mu_2)\in\mR^{d}\times{\cM_{\lambda^2}(\mR^d)}$
\ce
&&\langle{x_1-x_2,b(x_1,\mu_1)-b(x_2,\mu_2)}\rangle\\
&\leq&|x_1-x_2||b(x_1,\mu_1)-b(x_2,\mu_2)|\\
&\leq&\lambda_1\(|x_1-x_2|^2\gamma_1(|x_1-x_2|)+|x_1-x_2|\rho(\mu_1,\mu_2)\)\\
&\leq&\lambda_1\(|x_1-x_2|^2\gamma_1(|x_1-x_2|)+\frac{|x_1-x_2|^2}{2}+\frac{\rho^2(\mu_1,\mu_2)}{2}\)\\
&\leq&\lambda_1\(|x_1-x_2|^2\gamma_1(|x_1-x_2|)+|x_1-x_2|^2+\rho^2(\mu_1,\mu_2)\).
\de
Besides, by the proof of Theorem 2.3 in \cite{zx}, we know that there exists a $0<\eta<\frac{1}{e}$ such that
\ce
x^2\gamma_i(x)&\leq&\kappa_\eta(x^2), \qquad i=1,2, 
\de
where 
$$
\kappa_\eta(x):=\left\{
\begin{array}{rcl}
&0, &x=0,\\
&x\log x^{-1},             &{0< x\leq\eta},\\
&(\log \eta^{-1}-1)x+\eta,     &{x>\eta}, 
\end{array} \right.
$$
is a positive, strictly increasing, continuous concave function and satisfies $\kappa_\eta(0)=0$,  $\int_{0^+}\frac{1}{\kappa_\eta(x)+x}dx=\infty$. Thus, 
$$
\langle{x_1-x_2,b(x_1,\mu_1)-b(x_2,\mu_2)}\rangle\leq\lambda_1\(\kappa_\eta(|x_1-x_2|^2)+|x_1-x_2|^2+\rho^2(\mu_1,\mu_2)\).
$$
If $\sigma(x,\mu)$ satisfies ($\bf{H}_2'$), by the similar deduction to above it holds that 
\ce
\parallel{\sigma(x_1,\mu_1)-\sigma(x_2,\mu_2)}\parallel^2&\leq&\lambda_2\(|x_1-x_2|^2\gamma_2(|x_1-x_2|)+\rho^2(\mu_1,\mu_2)\)\\
&\leq&\lambda_2\(\kappa_\eta(|x_1-x_2|^2)+\rho^2(\mu_1,\mu_2)\).
\de
That is, ($\bf{H}_2'$) implies ($\bf{H}_2$). 
\er

\section{The existence and the uniqueness of strong solutions}\label{exun}

In the section, we study the existence and the uniqueness of strong solutions for Eq.$(\ref{eq1})$. The main result is the following theorem.

\bt\label{et0}
Suppose that $(\bf{H}_1)$-$(\bf{H}_2)$ hold and $\mE|\xi|^{2p}<\infty$ for any $p>1$. Then Eq.$(\ref{eq1})$ has a unique strong solution.
\et

The proof of the above theorem is made up of two parts--the existence  and the pathwise uniqueness of weak solutions. Firstly, we prove the existence of weak solutions for Eq.$(\ref{eq1})$. To do that, we introduce martingale solutions for Eq.$(\ref{eq1})$. Set
\ce
&&\cW:=C([0,T],\mR^d), \quad \sW=\sB(\cW),\\
&&\cW_t:=C([0,t],\mR^d), \quad \bar{\sW}_t=\cap_{s>t}\sB(\cW_t), \quad t\in[0,T].
\de

\bd\label{sd1}
A probability measure $P$ on $(\cW, \sW)$ is called a martingale solution of Eq.$(\ref{eq1})$ with the initial law $\mu_0$, if
\be
M_t^f&:=&f(w_t)-f(w_0)-\int_0^t(\mathscr{A}(\mu_s)f)(w_s)ds\label{eq2}, \quad\qquad f\in{C_0^2(\mR^d)},
\ee
is a continuous $\bar{\sW}_t$-adapted martingale, where $\mu_s:=P\circ{w_s^{-1}}$ denotes the law of $w_s$ under $P$ and
\ce
(\mathscr{A}(\mu){f})(x)&:=&\frac{1}{2}(\sigma(x,\mu)\sigma^*(x,\mu))^{ij}\partial_{ij}^2f+{b^i(x,\mu)\partial_i}f.
\de
\ed
Here and hereafter we use the convention that the repeated indices stand for the summation. We have the relationship between martingale solutions and weak solutions as follows.

\bp\label{er}
The existence of martingale solutions implies the existence of weak solutions and viceversa.
\ep

Since its proof is similar to that of \cite[Proposition 2.10]{chm}, we omit it. Next, we give a lemma which will take an important part in the sequel.

\bl\label{1}
Suppose $b(x,\mu)$ and $\sigma(x,\mu)$ satisfy $(\bf{H}_1)$.
If ($\hat{\cS}; \hat{W}, \hat{X}$) is a weak solution to Eq.$(\ref{eq1})$, where $\hat{\mE}(\cdot):=\mE^{\hat{\mP}}{(\cdot)}$ denotes the expectation under $\hat{\mP}$, it follows that for $p\geq1$
\be
\hat{\mE}({\mid{\hat{X}_t}\mid^{2p}})&\leq&{C(1}+\hat{\mE}\mid{\hat{X}_0}\mid^{2p})e^{Ct},\quad \qquad\qquad\qquad 0\leq{t}\leq{T},\label{2}\\
\hat{\mE}(\mid{\hat{X}_t-\hat{X}_s}\mid^{2p})&\leq&{C(1}+\hat{\mE}\mid{\hat{X}_0}\mid^{2p})(t-s)^p,\qquad\qquad 0\leq{s}<{t}\leq{T},\label{5}
\ee
where $C>0$ is a constant depending on $T$, $p$, $L_1$.
\el
\begin{proof}
Set $\tau_k:=inf\{t\geq0,\mid{\hat{X}_t}\mid\geq{k}\}$, $k\in\mN$. If these inequalities $(\ref{2})$ and $(\ref{5})$ hold for the process $\hat{X}_{\tau_k}$, let $k\rightarrow+\infty$, by Fatou's Lemma it follows that these inequalities $(\ref{2})$ and $(\ref{5})$ also hold for $\hat{X}_t$. So we might as well suppose that $\hat{X}_t$ is bounded.

For Eq.$(\ref{eq1})$, by the H\"{o}lder inequality and BDG inequality, it holds that
\ce
\hat{\mE}\mid{\hat{X}_t}\mid^{2p}&\leq&{3^{2p-1}}\(\hat{\mE}\mid{\hat{X}_0}\mid^{2p}+\hat{\mE}\mid\int_0^t{b(\hat{X}_s,\hat{\mu}_s)ds}\mid^{2p}+\hat{\mE}\mid\int_0^t{\sigma(\hat{X}_s,\hat{\mu}_s)d\hat{W}_s}\mid^{2p}\)\\
&\leq&{3^{2p-1}}\(\hat{\mE}\mid{\hat{X}_0}\mid^{2p}+\hat{\mE}\mid\int_0^t{b(\hat{X}_s,\hat{\mu}_s)ds}\mid^{2p}\\
&&+[p(2p-1)]^p\hat{\mE}\(\int_0^t\|\sigma(\hat{X}_s,\hat{\mu}_s)\|^2ds\)^p\)\\
&\leq&{3^{2p-1}}\(\hat{\mE}\mid{\hat{X}_0}\mid^{2p}+t^{2p-1}\(\int_0^t\hat{\mE}\mid b(\hat{X}_s,\hat{\mu}_s)\mid^{2p}ds\)\\
&&+[p(2p-1)]^pt^{p-1}\(\int_0^t\hat{\mE}\parallel\sigma(\hat{X}_s,\hat{\mu}_s)\parallel^{2p}ds\)\)\\
&\leq&C\(\hat{\mE}\mid{\hat{X}_0}\mid^{2p}+\int_0^t\hat{\mE}(\mid{b(\hat{X}_s,\hat{\mu}_s)}\mid^{2p}+\parallel{\sigma(\hat{X}_s,\hat{\mu}_s)}\parallel^{2p})ds\)\\
&\leq&C\(\hat{\mE}\mid{\hat{X}_0}\mid^{2p}+\int_0^t\hat{\mE}(1+\mid{\hat{X}_s}\mid^{2p}+\parallel{\hat{\mu}_s}\parallel^{2p}_{\lambda^2})ds\)\\
&\leq&C\(\hat{\mE}\mid{\hat{X}_0}\mid^{2p}+\int_0^t\hat{\mE}(1+\mid{\hat{X}_s}\mid^{2p}+\hat{\mE}(1+\mid{\hat{X}_s}\mid^{2p}))ds\)\\
&\leq&C\left(1+\hat{\mE}\mid{\hat{X}_0}\mid^{2p}+\int_0^t\hat{\mE}\mid \hat{X}_s\mid^{2p}ds\right),\quad\quad\quad 0\leq t\leq T,
\de
where $C>0$ is a constant depending on $T$, $p$, $L_1$. By Gronwall's inequality, one can get $(\ref{2})$.

By the similar deduction to above, we obtain
\ce
\hat{\mE}\mid \hat{X}_t-\hat{X}_s\mid^{2p}&\leq& C_p\hat{\mE}\(\mid\int_s^tb(\hat{X}_u,\hat{\mu}_u)du\mid^{2p}+\mid\int_s^t\sigma(\hat{X}_u,\hat{\mu}_u)d\hat{W}_u\mid^{2p}\)\\
&\leq&C_{p,T}(t-s)^{p-1}\int_s^t\hat{\mE}(\mid{b(\hat{X}_u,\hat{\mu}_u)}\mid^{2p}+\|\sigma(\hat{X}_u,\hat{\mu}_u)\|^{2p})du\\
&\leq&C_{p,T}(t-s)^{p-1}\int_s^t(1+\hat{\mE}\mid \hat{X}_u\mid^{2p})du\\
&\leq&C(1+\hat{\mE}\mid{\hat{X}_0}\mid^{2p})(t-s)^p,\quad\quad\quad 0\leq s< t\leq T.
\de
The proof is completed.
\end{proof}

\bp\label{et}
Suppose that $(\bf{H}_1)$ holds and $\mE|\xi|^{2p}<\infty$ for any $p>1$. Then there exists a martingale solution to Eq.$(\ref{eq1})$.
\ep
\begin{proof}
Firstly, for fixed ${n\in{\mN}}$, consider the following Euler-Maruyama approximation equation
\be
d{X_t^n}=b(X_{t_n}^n,\mu_{t_n}^n)dt+\sigma(X_{t_n}^n,\mu_{t_n}^n)d{W_t},\label{aeq2}
\ee
where ${X_0^n=\xi},t_n=\frac{[2^nt]}{2^n}$ and $[a]$ denotes the integer part of $a$.
By solving a deterministic problem, this equation can be solved step by step. That is, there exists a solution $X^n$ to Eq.$(\ref{aeq2})$. By $(\ref{eq3})$ and Lemma $\ref{1}$, we have
\be
\mE({\mid{X_{t}^n}\mid^{2p}})&\leq&{C(1}+\mE\mid{\xi}\mid^{2p})e^{Ct},\quad\qquad\qquad 0\leq{t}\leq{T}\label{4},\\
\mE(\mid{X_t^n-X_s^n}\mid^{2p})&\leq&{C(1}+\mE\mid{\xi}\mid^{2p})(t-s)^p,\qquad 0\leq{s}<{t}\leq{T},\no
\ee
where $C$ is independent of $n$. Since $\mE\mid{\xi}\mid^{2p}<+\infty$, we further have
\ce
\sup_{n\geq1}\mE\mid{X_0^n}\mid^{2p}&=&{\mE\mid{\xi}\mid^{2p}}<+\infty,\\
\sup_{n\geq1}\mE(\mid{X_t^n-X_s^n}\mid^{2p})&\leq&{C(1}+\mE\mid{\xi}\mid^{2p})(t-s)^p\leq{C(t-s)^{p}}.
\de
Set $P^n:={\mP\circ(X^n)^{-1}}$, and then by Lemma 20.3 in  \cite[P.185]{hzy}  we derive that $\{P^n\}$ is tight. So there exist a subsequence still denoted by $\{P^n\}$ and $P^0$ such that $P^n$ weakly converges to $P^0$ as $n\rightarrow{+}\infty$.

Now set
\ce
M_t^{n,f}&:=&f(w_t)-f(w_0)
-\frac{1}{2}{\int_0^t{(\sigma(w_{s_n},\mu_{s_n}^n)\sigma^*(w_{s_n},\mu_{s_n}^n))^{ij}\partial_{ij}^2f(w_s)}ds}\\
&&-{\int_0^t{b^i(w_{s_n},\mu_{s_n}^n)\partial_i}f(w_s)ds},\quad\quad\quad\quad f\in{C_0^2(\mR^d)}.
\de
Since Eq.$(\ref{aeq2})$ has a weak solution $X^n$, by Proposition $\ref{er}$, we know that there exists a martingale solution $P^n$ on $(\cW,\sW)$ of Eq.$(\ref{aeq2})$, which yields that $M_t^{n,f}$ is a continuous $\bar{\sW}_t$-adapted martingale under $P^n$.  So for any continuous, bounded and $\bar{\sW}_{s}$-measurable functional $G_s$,
\ce
\mE^{P^n}((M_t^{n,f}-M_s^{n,f})G_s)=0,\quad\quad 0\leq{s}<{t}\leq{T}.
\de
To prove that $P^0$ on $(\cW,\sW)$ is a martingale solution to Eq.$(\ref{eq1})$, we just need to prove that $M_t^f$ defined by $(\ref{eq2})$ is a continuous $\bar{\sW}_t$-adapted martingale under $P^0$. That is,
\ce
\mE^{P^0}\((M_t^f-M_s^f)G_s\)=\int_{\cW}\left(\left(f(w_t)-f(w_s)-\int_s^t{\mathscr{A}(\mu_u)f(w_u)}du\right)G_s(w)\right)P^0(\dif w)=0.
\de

Note that $P^n$ weakly converges to $P^0$. Thus, it is clear that
\ce
\lim_{n\rightarrow\infty}\int_{\cW}\((f(w_t)-f(w_s))G_s(w)\)P^n(\dif w)=\int_{\cW}\((f(w_t)-f(w_s))G_s(w)\)P^0(\dif w).
\de
We now prove that
\be
&&\lim_{n\rightarrow\infty}\int_{\cW}\left(\left(\int_s^t{b^i(w_{u_n},\mu_{u_n}^n)\partial_i}f(w_u)du\right)G_s(w)\right)P^n(\dif w)\no\\
&=&\int_{\cW}\left(\left(\int_s^t{b^i(w_{u},\mu_{u})\partial_i}f(w_u)du\right)G_s(w)\right)P^0(\dif w),
\label{eq5}
\ee
and
\be
&&\lim_{n\rightarrow\infty}\int_{\cW}\left(\left(\int_s^t{(\sigma(w_{u_n},\mu_{u_n}^n)\sigma^*(w_{u_n},\mu_{u_n}^n))^{ij}\partial_{ij}^2f(w_u)}du\right)G_s(w)\right)P^n(\dif w)\no\\
&=&\int_{\cW}\left(\left(\int_s^t{(\sigma(w_{u},\mu_{u})\sigma^*(w_{u},\mu_{u}))^{ij}\partial_{ij}^2f(w_u)}du\right)G_s(w)\right)P^0(\dif w).
\label{eq6}
\ee
With the help of Theorem c.6 \cite [P.324]{hzy} and the weak convergence of $P^n$ to $P^0$, we know that there exist a probability space $(\tilde{\Omega},\tilde{\mathscr{F}},\tilde{\mP})$ and $\cW$-valued processes $\tilde{X}_{\cdot}^n$, $\tilde{X}_{\cdot}$ on it satisfying

(i) The law of $\tilde{X}_{\cdot}^n$ and $\tilde{X}_{\cdot}$ are $P^n$ and $P^0$, respectively,

(ii) $ \tilde{X}_{\cdot}^n\stackrel{a.s.}\rightarrow\tilde{X}_{\cdot}$ as $n\rightarrow\infty$.\\
Based on (i), (\ref{eq5}) (\ref{eq6}) become
\be
&&\lim_{n\rightarrow\infty}\mE^{\tilde{\mP}}\left(\left(\int_s^t{b^i(\tilde{X}^n_{u_n},\mu_{u_n}^n)\partial_i}f(\tilde{X}^n_u)du\right)G_s(\tilde{X}^n_{\cdot})\right)\no\\
&=&\mE^{\tilde{\mP}}\left(\left(\int_s^t{b^i(\tilde{X}_{u},\mu_{u})\partial_i}f(\tilde{X}_u)du\right)G_s(\tilde{X}_{\cdot})\right),
\label{eq7}
\ee
and
\be
&&\lim_{n\rightarrow\infty}\mE^{\tilde{\mP}}\left(\left(\int_s^t(\sigma(\tilde{X}^n_{u_n},\mu_{u_n}^n)\sigma^*(\tilde{X}^n_{u_n},\mu_{u_n}^n))^{ij}\partial_{ij}^2f(\tilde{X}^n_u)du\right)G_s(\tilde{X}^n_{\cdot})\right)\no\\
&=&\mE^{\tilde{\mP}}\left(\left(\int_s^t(\sigma(\tilde{X}_{u},\mu_{u})\sigma^*(\tilde{X}_{u},\mu_{u}))^{ij}\partial_{ij}^2f(\tilde{X}_u)du\right)G_s(\tilde{X}_{\cdot})\right).
\label{eq8}
\ee

In the following, we are devoted to proving (\ref{eq7}). On one side, by (ii), it holds that $\tilde{X}_{u_n}^n\stackrel{a.s.}\rightarrow\tilde{X}_u$ for $u\in[s,t]$ as $n\rightarrow\infty$. Next, we observe $\rho(\mu_{u_n}^n,\mu_{u})$. By the definition of $\rho$, it holds that
\ce
\rho(\mu_{u_n}^n,\mu_{u})&=&\sup_{\parallel{\varphi}\parallel_{C_\rho(\mR^d)\leq1}}\left|{\int_{\mR^d}\varphi(x)\mu_{u_n}^n(dx)-\int_{\mR^d}\varphi(x)\mu_{u}(dx)}\right|\\
&=&\sup_{\parallel{\varphi}\parallel_{C_\rho(\mR^d)\leq1}}\left|\mE^{\tilde{\mP}}\varphi(\tilde{X}_{u_n}^n)-\mE^{\tilde{\mP}}\varphi(\tilde{X}_{u})\right|\\
&\leq&\sup_{\parallel{\varphi}\parallel_{C_\rho(\mR^d)\leq1}}\mE^{\tilde{\mP}}\left|\varphi(\tilde{X}_{u_n}^n)-\varphi(\tilde{X}_{u})\right|\\
&\leq&\mE^{\tilde{\mP}}\left|\tilde{X}_{u_n}^n-\tilde{X}_{u}\right|.
\de
Note that for any $\lambda>0$,
\ce
\int_{|\tilde{X}_{u_n}^n|>\lambda}|\tilde{X}_{u_n}^n| d\tilde{\mP}=\int_{|\tilde{X}_{u_n}^n|>\lambda}\frac{|\tilde{X}_{u_n}^n|^p}{\lambda^p}\lambda^p d\tilde{\mP}\leq\int_{|\tilde{X}_{u_n}^n|>\lambda}\frac{|\tilde{X}_{u_n}^n|^{2p}}{\lambda^p} d\tilde{\mP}\leq\frac{1}{\lambda^p} \mE^{\tilde{\mP}}\left|\tilde{X}_{u_n}^n\right|^{2p}=\frac{1}{\lambda^p}\mE\left|X_{u_n}^n\right|^{2p}.
\de
Thus, by $(\ref{4})$ we have that 
\ce
\lim_{\lambda\rightarrow\infty}\sup_{n\geq 1}\int_{|\tilde{X}_{u_n}^n|>\lambda}|\tilde{X}_{u_n}^n| d\tilde{\mP}=0,
\de
and then $\{\tilde{X}_{u_n}^n, n\geq 1\}$ is uniformly integrable. Based on \cite[Theorem 4.5]{hzy}, one can know that uniform integrability of $\{\tilde{X}_{u_n}^n, n\geq 1\}$ and almost sure convergence of $\tilde{X}_{u_n}$ to $\tilde{X}_u$ imply that $\lim_{n\rightarrow\infty}\mE^{\tilde{\mP}}\left|\tilde{X}_{u_n}^n-\tilde{X}_{u}\right|=0$ and furthermore $\lim_{n\rightarrow\infty}\rho(\mu_{u_n}^n,\mu_{u})=0$.

On the other side, by (i) $(\ref{eq3})$ and $(\ref{4})$, it holds that
\ce
\mE^{\tilde{P}}\mid{b(\tilde{X}^n_{u_n},\mu_{u_n}^n)}\mid&=&\mE\mid{b(X^n_{u_n},\mu^n_{u_n})}\mid\leq\(\mE\mid{b(X^n_{u_n},\mu^n_{u_n})}\mid^{2p}\)^{1/2p}\\
&\leq&C\(\mE\(1+\mid{X^n_{u_n}}\mid^{2p}+\parallel{\mu^n_{u_n}}\parallel^{2p}_{\lambda^2}\)\)^{1/2p}\\
&\leq&C\(\mE\(1+\mid{X^n_{u_n}}\mid^{2p}+\mE(1+\mid{X^n_{u_n}}\mid^{2p})\)\)^{1/2p}\\
&\leq&C\(\mE(1+\mid{X^n_{u_n}}\mid^{2p})\)^{1/2p}<\infty.
\de
By the continuity of $b$ and the dominated convergence theorem, we obtain $(\ref{eq7})$. By the similar means, one can prove $(\ref{eq8})$. The proof is now completed.
\end{proof}

So, by Proposition $\ref{er}$ and Theorem $\ref{et}$, we know that Eq.$(\ref{eq1})$ has a weak solution. Next, we prove that pathwise uniqueness holds for Eq.$(\ref{eq1})$ under certain conditions. The following lemma is known (c.f. \cite[Lemma 116, P.79 and Lemma 144, P.113]{SR}). For the readers' convenience, we give a short proof.

\bl\label{3}
For any $t\geq0$, if $y_t$ satisfies $0\leq{y_t}\leq\int_0^t\left(\kappa_1(y_s)+\kappa_2(y_s)\right)ds<\infty$, where $\kappa(u)$ satisfies the conditions in $(\bf{H}_2)$, then $y_t\equiv0,\forall{t\geq0}$.
\el
\begin{proof}
Set $z_t:=\int_0^t(\kappa_1(y_s)+\kappa_2(y_s))ds$, and then we just need to prove $z_t=0$. Note that $z_t$ is absolutely continuous and nondecreasing. Thus, it holds that
\ce
\frac{d{z_t}}{dt}=\kappa_1(y_t)+\kappa_2(y_t)\leq{\kappa_1(z_t)+\kappa_2(z_t)}.
\de
Let $t_0:=\sup\{{t\geq0;z_s=0,\forall{s\in[0,t]}}\}$. If $t_0<\infty$, then $z_t>0$, $t>t_0$.  Therefore, we have
\ce
\infty=\int_0^{z_{t_0+\varepsilon}}\frac{du}{\kappa_1(u)+\kappa_2(u)}=\int_{t_0}^{t_0+\varepsilon}\frac{dz_t}{\kappa_1(z_t)+\kappa_2(z_t)}\leq\int_{t_0}^{t_0+\varepsilon}dt\leq\varepsilon,\quad\quad\forall\varepsilon>0,
\de
which is a contradiction. So $t_0=\infty$ and $z_t=0$.
\end{proof}

\bp\label{et2}
Suppose that $(\bf{H}_2)$ holds. Then the pathwise uniqueness holds for Eq.$(\ref{eq1})$.
\ep
\begin{proof}
Suppose that ($\hat{\cS}; \hat{W}, \hat{X}^1$) and ($\hat{\cS}; \hat{W}, \hat{X}^2$) are two weak solutions to Eq.$(\ref{eq1})$ with $\hat{X}^1_0=\hat{X}^2_0$. Set
\ce
Z_t&:= &\hat{X}_t^1-\hat{X}_t^2,
\de
and then $Z_t$ satisfies
\ce
Z_t=\int_0^t\(b(\hat{X}_s^1,\hat{\mu}_s^1)-b(\hat{X}_s^2,\hat{\mu}_s^2)\)ds+\int_0^t\(\sigma(\hat{X}_s^1,\hat{\mu}_s^1)-\sigma(\hat{X}_s^2,\hat{\mu}_s^2)\)d\hat{W}_s.
\de
Applying the It\^o formula to $|Z_t|^2$, we obtain that
\ce
|Z_t|^2&=&\int_0^t2\langle{Z_s,b(\hat{X}_s^1,\hat{\mu}_s^1)-b(\hat{X}_s^2,\hat{\mu}_s^2})\rangle ds+\int_0^t\parallel{\sigma(\hat{X}_s^1,\hat{\mu}_s^1)-\sigma(\hat{X}_s^2,\hat{\mu}_s^2)}\parallel^2ds\\
&+&\int_0^t2\langle{Z_s,(\sigma(\hat{X}_s^1,\hat{\mu}_s^1)-\sigma(\hat{X}_s^2,\hat{\mu}_s^2))d\hat{W}_s}\rangle.
\de
By taking the expectation on two sides, one can have
\ce
\hat{\mE}|Z_t|^2=\hat{\mE}\int_0^t2\langle{Z_s,b(\hat{X}_s^1,\hat{\mu}_s^1)-b(\hat{X}_s^2,\hat{\mu}_s^2})\rangle ds+\hat{\mE}\int_0^t\parallel{\sigma(\hat{X}_s^1,\hat{\mu}_s^1)-\sigma(\hat{X}_s^2,\hat{\mu}_s^2)}\parallel^2ds.
\de
Put $G_t:=\hat{\mE}\mid{Z_t}\mid^2$, and by $(\bf{H}_2)$, it holds that
\ce
G_t&=&\hat{\mE}\int_0^t\(2\langle{Z_s,b(\hat{X}_s^1,\hat{\mu}_s^1)-b(\hat{X}_s^2,\hat{\mu}_s^2})\rangle +\parallel{\sigma(\hat{X}_s^1,\hat{\mu}_s^1)-\sigma(\hat{X}_s^2,\hat{\mu}_s^2)}\parallel^2\)ds\\
&\leq&L_2\hat{\mE}\int_0^t\(\kappa_1(|Z_s|^2)+\kappa_2\left(\rho^2(\hat{\mu}_s^1,\hat{\mu}_s^2)\right)\)ds.
\de
Note that
\be
\rho(\hat{\mu}_s^1,\hat{\mu}_s^2)&=&\sup_{\parallel{\varphi}\parallel_{C_\rho(\mR^d)\leq1}}\left|{\int_{\mR^d}\varphi(x)\hat{\mu}_s^1(dx)-\int_{\mR^d}\varphi(x)\hat{\mu}_s^2(dx)}\right|\no\\
&=&\sup_{\parallel{\varphi}\parallel_{C_\rho(\mR^d)\leq1}}\left|\hat{\mE}\varphi(\hat{X}_s^1)-\hat{\mE}\varphi(\hat{X}_s^2)\right|\no\\
&\leq&\sup_{\parallel{\varphi}\parallel_{C_\rho(\mR^d)\leq1}}\hat{\mE}\left|\varphi(\hat{X}_s^1)-\varphi(\hat{X}_s^2)\right|\no\\
&\leq&\hat{\mE}\left|\hat{X}_s^1-\hat{X}_s^2\right|,\label{8}
\ee
and
$$
\rho^2(\hat{\mu}_s^1,\hat{\mu}_s^2)\leq\(\hat{\mE}\left|\hat{X}_s^1-\hat{X}_s^2\right|\)^2\leq\hat{\mE}\left|\hat{X}_s^1-\hat{X}_s^2\right|^2=\hat{\mE}|Z_s|^2=G_s.
$$
Thus, by the Jensen inequality, we get that
\ce
G_t&\leq& L_2\hat{\mE}\int_0^t\(\kappa_1(|Z_s|^2)+\kappa_2(G_s)\)ds
\leq L_2\int_0^t\(\kappa_1(\hat{\mE}|Z_s|^2)+\kappa_2(G_s)\)ds\\
&=&L_2\int_0^t\(\kappa_1(G_s)+\kappa_2(G_s)\)ds.
\de
By Lemma $\ref{3}$ we have that $G_t=0$ and then $Z_t=0,\forall{t\geq0}$, a.s.. Therefore the pathwise uniqueness is right.
\end{proof}

Finally, Theorem $\ref{et0}$ can be proved by Theorem $\ref{et}$, Theorem $\ref{et2}$ and \cite[Proposition 3.20, P.309]{dz}.

\section{The convergence rate for the Euler-Maruyama approximation}\label{cora}

In the section we consider the convergence rate for the Euler-Maruyama approximation $\{X_t^n\}$ defined in $(\ref{aeq2})$, i.e.
\ce
X_t^n=\xi+\int_0^tb(X_{s_n}^n,\mu_{s_n}^n)ds+\int_0^t\sigma(X_{s_n}^n,\mu_{s_n}^n)d{W_s}.
\de
where $s_n=\frac{[2^n s]}{2^n}$ and $[a]$ denotes the integer part of $a$.

\bt\label{et3}
Suppose $b$ and $\sigma$ satisfy $(\bf{H}_1)$ and $(\bf{H}'_2)$ and $\mE|\xi|^{2p}<\infty$ for any $p>1$. Then there exists a $T_0>0$ such that
\ce
\mE\(\sup_{t\in[0,T]}\mid X_t^n-X_t\mid^2\)=O(2^{-n}T_0),
\de
where $O(2^{-n}T_0)$ means that $\frac{O(2^{-n}T_0)}{2^{-n}T_0}$ is bounded.
\et
\begin{proof}
Set $H_t:=X_t^n-X_t$, and then $H_t$ satisfies
\ce
H_t=\int_0^t\(b(X_{s_n}^n,\mu_{s_n}^n)-b(X_s,\mu_s)\)ds+\int_0^t\(\sigma(X_{s_n}^n,\mu_{s_n}^n)-\sigma(X_s,\mu_s)\)dW_s.
\de
It follows from the It\^{o} formula that
\ce
\mid H_t\mid^2=J_1+J_2+J_3,
\de
where
\ce
J_1&:=&\int_0^t2\langle H_s,b(X_{s_n}^n,\mu_{s_n}^n)-b(X_s,\mu_s)\rangle ds,\\
J_2&:=&\int_0^t2\langle H_s,(\sigma(X_{s_n}^n,\mu_{s_n}^n)-\sigma(X_s,\mu_s))dW_s\rangle,\\
J_3&:=&\int_0^t\parallel\sigma(X_{s_n}^n,\mu_{s_n}^n)-\sigma(X_s,\mu_s)\parallel^2ds.
\de

For $J_1$, by $(\bf{H}'_2)$ and $(\ref{8})$ it holds that 
\ce
\mE|J_1|&\leq&2\mE\int_0^t\mid H_s\mid\mid b(X_{s_n}^n,\mu_{s_n}^n)-b(X_s,\mu_s)\mid ds\\
&\leq&2\mE\int_0^t\(\mid H_s\mid\mid b(X_{s_n}^n,\mu_{s_n}^n)-b(X_{s}^n,\mu_{s}^n)\mid +\mid H_s\mid\mid b(X_{s}^n,\mu_{s}^n)-b(X_s,\mu_s)\mid \)ds\\
&\leq&\mE\int_0^t\(\mid H_s\mid^2+\mid b(X_{s_n}^n,\mu_{s_n}^n)-b(X_s^n,\mu_s^n)\mid^2\) ds\\
&&+2\mE\int_0^t\lambda_1\(\mid H_s\mid^2\gamma_1(\mid H_s\mid)+\mid H_s\mid\rho(\mu_s^n,\mu_s)\)ds\\
&\leq&\mE\int_0^t\mid H_s\mid^2 ds+2\mE\int_0^t\lambda_1^2\(\mid X_{s_n}^n-X_s^n\mid^2\gamma_1^2(\mid X_{s_n}^n-X_s^n\mid)+\rho^2(\mu_{s_n}^n,\mu_s^n)\)ds\\
&&+2\mE\int_0^t\lambda_1\mid H_s\mid^2\gamma_1(\mid H_s\mid)ds+\lambda_1\mE\int_0^t\mid H_s\mid^2ds+\lambda_1\mE\int_0^t\rho^2(\mu_s^n,\mu_s)ds\\
&\leq&C\mE\int_0^t\mid H_s\mid^2 ds+C_{\lambda_1}\mE\int_0^t\kappa_\eta^2(\mid X_{s_n}^n-X_s^n\mid)ds+C_{\lambda_1}\int_0^t\mE(\mid X_{s_n}^n-X_s^n\mid^2)ds\\
&&+C_{\lambda_1}\mE\int_0^t\kappa_\eta(\mid H_s\mid^2)ds,
\de
where in the last inequality the following result is used that
\ce
x\gamma_1(x)&\leq&\kappa_\eta(x), \qquad x>0,\\
x^2\gamma_1(x)&\leq&\kappa_\eta(x^2),
\de
and for $0<\eta<\frac{1}{e}$
$$
\kappa_\eta(x)=\left\{
\begin{array}{rcl}
&0, &x=0,\\
&x\log x^{-1},             &{0< x\leq\eta},\\
&(\log \eta^{-1}-1)x+\eta,     &{x>\eta}.
\end{array} \right.
$$
Here the properties of $\kappa_\eta$ can be referred to in Remark \ref{h2h2}. And then, the Jensen inequality gives that
\ce
\mE|J_1|&\leq&C\int_0^t\mE\mid H_s\mid^2 ds+C_{\lambda_1}\int_0^t\kappa_\eta^2((\mE\mid X_{s_n}^n-X_s^n\mid^2)^{1/2})ds\no\\
&&+C_{\lambda_1}\int_0^t\mE(\mid X_{s_n}^n-X_s^n\mid^2)ds+C_{\lambda_1}\int_0^t\kappa_\eta(\mE\mid H_s\mid^2)ds,
\de
and furthermore
\be
\mE\left(\sup_{t\in[0,T]}|J_1|\right)&\leq&C\int_0^T\mE\left(\sup_{r\in[0,s]}\mid H_r\mid^2\right) ds+C_{\lambda_1}\int_0^T\kappa_\eta^2\left(\left(\mE\left(\sup_{r\in[0,s]}\mid X_{r_n}^n-X_r^n\mid^2\right)\right)^{1/2}\right)ds\no\\
&&+C_{\lambda_1}\int_0^T\mE\left(\sup_{r\in[0,s]}\mid X_{r_n}^n-X_r^n\mid^2\right)ds+C_{\lambda_1}\int_0^T\kappa_\eta\left(\mE\left(\sup_{r\in[0,s]}\mid H_r\mid^2\right)\right)ds.\no\\
\label{6}
\ee
Similarly, we obtain
\be
\mE\left(\sup_{t\in[0,T]}|J_3|\right)&\leq&C\int_0^T\mE\left(\sup_{r\in[0,s]}\mid H_r\mid^2\right) ds+C_{\lambda_2}\int_0^T\kappa_\eta\left(\mE\left(\sup_{r\in[0,s]}\mid X_{r_n}^n-X_r^n\mid^2\right)\right)ds\no\\
&&+C_{\lambda_2}\int_0^T\mE\left(\sup_{r\in[0,s]}\mid X_{r_n}^n-X_r^n\mid^2\right)ds+C_{\lambda_2}\int_0^T\kappa_\eta\left(\mE\left(\sup_{r\in[0,s]}\mid H_r\mid^2\right)\right)ds.\no\\
\label{7}
\ee
For $J_2$, by $(\bf{H}'_2)$, $(\ref{7})$, the BDG inequality and the Young inequality, one can get that
\be
\mE\left(\sup_{t\in[0,T]}|J_2|\right)&\leq& C\mE\(\int_0^T\mid H_s\mid^2\parallel\sigma(X_{s_n}^n,\mu_{s_n}^n)-\sigma(X_s,\mu_s)\parallel^2ds\)^\frac{1}{2}\no\\
&\leq&C\mE\(\sup_{t\in[0,T]}\mid H_t\mid^2\int_0^T\parallel\sigma(X_{s_n}^n,\mu_{s_n}^n)-\sigma(X_s,\mu_s)\parallel^2ds\)^\frac{1}{2}\no\\
&\leq&\frac{1}{4}\mE\left(\sup_{t\in[0,T]}\mid H_t\mid^2\right)+C\mE\int_0^T\parallel\sigma(X_{s_n}^n,\mu_{s_n}^n)-\sigma(X_s,\mu_s)\parallel^2ds\no\\
&\leq&\frac{1}{4}\mE\left(\sup_{t\in[0,T]}\mid H_t\mid^2\right)+C\int_0^T\mE\left(\sup_{r\in[0,s]}\mid H_r\mid^2\right) ds\no\\
&&+C_{\lambda_2}\int_0^T\kappa_\eta\left(\mE\left(\sup_{r\in[0,s]}\mid X_{r_n}^n-X_r^n\mid^2\right)\right)ds\no\\
&&+C_{\lambda_2}\int_0^T\mE\left(\sup_{r\in[0,s]}\mid X_{r_n}^n-X_r^n\mid^2\right)ds\no\\
&&+C_{\lambda_2}\int_0^T\kappa_\eta\left(\mE\left(\sup_{r\in[0,s]}\mid H_r\mid^2\right)\right)ds.
\label{11}
\ee

Combining $(\ref{6})$-$(\ref{11})$, we know that
\ce
\mE\left(\sup_{t\in[0,T]}\mid H_t\mid^2\right)&\leq&C\int_0^T\mE\left(\sup_{r\in[0,s]}\mid H_r\mid^2\right) ds+C\int_0^T\kappa^2_\eta\left(\left(\mE\left(\sup_{r\in[0,s]}\mid X_{r_n}^n-X_r^n\mid^2\right)\right)^{\frac{1}{2}}\right)ds\\
&&+C\int_0^T\kappa_\eta\left(\mE\left(\sup_{r\in[0,s]}\mid X_{r_n}^n-X_r^n\mid^2\right)\right)ds+C\int_0^T\mE\left(\sup_{r\in[0,s]}\mid X_{r_n}^n-X_r^n\mid^2\right)ds\\
&&+C\int_0^T\kappa_\eta\left(\mE\left(\sup_{r\in[0,s]}\mid H_r\mid^2\right)\right)ds,
\de
where $C>0$ is a constant depending on $\lambda_1$, $\lambda_2$. Next, we estimate $\mE\left(\sup\limits_{r\in[0,s]}\mid X_{r_n}^n-X_r^n\mid^2\right)$. Note that for $r_n=\frac{i}{2^n}T\leq r<\frac{i+1}{2^n}T, i=0,1,2,\cdots,2^n-1$
\ce
X_r^n&=&X^n_{r_n}+\int_{r_n}^rb(X_{s_n}^n,\mu_{s_n}^n)ds+\int_{r_n}^r\sigma(X_{s_n}^n,\mu_{s_n}^n)d{W_s}\\&=&X^n_{r_n}+b(X_{r_n}^n,\mu_{r_n}^n)(r-r_n)+\sigma(X_{r_n}^n,\mu_{r_n}^n)(W_r-W_{r_n}).
\de
By $(\ref{eq3})$, $(\ref{4})$, the H\"older inequality and the BDG inequality, it holds that
\ce
&&\mE\left(\sup_{\frac{i}{2^n}T\leq r<\frac{i+1}{2^n}T}|X_r^n-X^n_{r_n}|^2\right)\\
&\leq&2\mE\left(|b(X_{r_n}^n,\mu_{r_n}^n)|^2\left|\frac{i+1}{2^n}T-\frac{i}{2^n}T\right|^2\right)\\
&&+2\mE\left(\|\sigma(X_{r_n}^n,\mu_{r_n}^n)\|^2\sup_{\frac{i}{2^n}T\leq r<\frac{i+1}{2^n}T}|W_r-W_{r_n}|^2\right)\\
&\leq&C2^{-n}T+C\left(\mE\|\sigma(X_{r_n}^n,\mu_{r_n}^n)\|^4\right)^{1/2}\left(\mE\sup_{\frac{i}{2^n}T\leq r<\frac{i+1}{2^n}T}|W_r-W_{r_n}|^4\right)^{1/2}\\
&\leq&C2^{-n}T,
\de
where the constant $C>0$ is independent of $n$. Thus, we obtain
\ce
\mE\left(\sup_{t\in[0,T]}\mid H_t\mid^2\right)&\leq&C\int_0^T\mE\left(\sup_{r\in[0,s]}\mid H_r\mid^2\right) ds+C\int_0^T\kappa_\eta\left(\mE\left(\sup_{r\in[0,s]}\mid H_r\mid^2\right)\right)ds\\
&&+CT\kappa_\eta^2(C(2^{-n}T)^{1/2})+CT\kappa_\eta(C2^{-n}T)+CT(2^{-n}T)\\
&\leq&C\int_0^T\left(\mE\left(\sup_{r\in[0,s]}\mid H_r\mid^2\right)+\kappa_\eta\left(\mE\left(\sup_{r\in[0,s]}\mid H_r\mid^2\right)\right)\right)ds\\
&&+CT\kappa_\eta^2(C(2^{-n}T)^{1/2})+CT\kappa_\eta(C2^{-n}T)+CT(2^{-n}T).
\de
By Lemma 144 in \cite[P.113]{SR} and Lemma 2.1 in \cite{zx} we have
\ce
\mE\left(\sup_{t\in[0,T]}\mid H_t\mid^2\right)\leq A^{\exp\{-CT\}},
\de
where $A:=CT\kappa_\eta^2(C(2^{-n}T)^{1/2})+CT\kappa_\eta(C2^{-n}T)+CT(2^{-n}T)$. Thus, there exists a $T_0>0$ such that
\ce
\mE\(\sup_{t\in[0,T_0]}\mid X_t^n-X_t\mid^2\)=O(2^{-n}T_0).
\de
If $T_0\geq T$, the proof is over; if $T_0<T$, on $[T_0, 2T_0], [2T_0, 3T_0], \cdots, [[\frac{T}{T_0}]T_0,T]$, by the same way to the above we deduce and conclude that
\ce
\mE\(\sup_{t\in[0,T]}\mid X_t^n-X_t\mid^2\)=O(2^{-n}T_0).
\de
The proof is completed.
\end{proof}

\bigskip

\textbf{Acknowledgements:}

Two authors would like to thank Professor Xicheng Zhang for his valuable discussions. And they would also wish to thank the anonymous referee for giving useful suggestions to improve this paper.

\end{document}